\documentclass[11pt]{article}
\usepackage[margin=1in]{geometry}
\usepackage{amsmath,amssymb,amsthm,mathtools}
\usepackage{microtype}
\usepackage[dvipsnames]{xcolor}
\usepackage{hyperref}
\usepackage{enumitem}
\usepackage[all]{nowidow}
\hypersetup{colorlinks=true,linkcolor=blue!50!black,citecolor=blue!50!black,urlcolor=blue!50!black}

\newtheorem{theorem}{Theorem}[section]
\newtheorem{proposition}[theorem]{Proposition}
\newtheorem{lemma}[theorem]{Lemma}
\newtheorem{corollary}[theorem]{Corollary}
\theoremstyle{definition}
\newtheorem{definition}[theorem]{Definition}
\newtheorem{example}[theorem]{Example}
\newtheorem{remark}[theorem]{Remark}

\makeatletter
\def\@fnsymbol#1{\ensuremath{\ifcase#1\or \mid\or \dagger\or \ddagger\or
   \mathsection\or \mathparagraph\or \|\or **\or \dagger\dagger
   \or \ddagger\ddagger \else\@ctrerr\fi}}
\makeatother

\newcommand{\g}{\mathfrak{g}}
\newcommand{\h}{\mathfrak{h}}
\newcommand{\K}{\mathbb{K}}
\newcommand{\id}{\mathrm{id}}
\newcommand{\End}{\mathrm{End}}
\newcommand{\lr}[1]{\langle #1\rangle}

\newcommand{\nc}{\newcommand}
\nc{\yuan}[1]{\textcolor{orange}{yuan:#1}}

\title{Admissible Yang--Baxter equation for Nijenhuis perm superalgebras}

\author{
Imed Basdouri\thanks{Department of Mathematics, Faculty of Sciences,
University of Gafsa, Gafsa, Tunisia.
\texttt{basdourimed@yahoo.fr}}
\and Bouzid Mosbahi\thanks{Department of Mathematics, Faculty of Sciences,
University of Sfax, Sfax, Tunisia.
\texttt{mosbahi.bouzid.etud@fss.usf.tn}}
\and Yuanyuan Zhang\textsuperscript{$\ddagger*$}
}
\footnotetext{\textsuperscript{$*$} Corresponding author.}
\footnotetext{\textsuperscript{$\ddagger$} School of Mathematics and Statistics, Henan University, Henan, Kaifeng 475004, P. R. China. \texttt{zhangyy17@henu.edu.cn}}
\date{}

\allowdisplaybreaks
\begin{document}
\maketitle

\begin{abstract}
We develop a superization of the theory of admissible Yang--Baxter equations for Nijenhuis perm algebras initiated by Ma and Song. We introduce Nijenhuis perm superalgebras, their representations, and the dual representation, with Koszul signs chosen so that the contraction of the Yang--Baxter tensor matches the $\mathcal O$-operator equation exactly. We define Nijenhuis perm superbialgebras and characterise them by matched pairs. In the coboundary case, we show that every even supersymmetric solution of the $S$-admissible perm Yang--Baxter equation gives rise to a Nijenhuis perm superbialgebra, and we give sufficient conditions for a perm superalgebra to admit a Nijenhuis operator constructed from such a solution. Examples with purely odd components illustrate the role of the grading.
\end{abstract}


\section{Introduction}

Nijenhuis operators play a role in the deformation theory of algebras, in bialgebra theory, in integrable systems, in almost complex structures and in Nijenhuis geometry. A (right) perm algebra is an associative algebra satisfying
\[
(xy)z=x(yz)=x(zy).
\]
Perm algebras have been studied in connection with differential algebras and metabelian Lie algebras \cite{KS,MaSa}. Perm bialgebras, and the associated \emph{perm Yang--Baxter equation} (called the $S$-equation in \cite{Hou}), were introduced by Hou \cite{Hou} through Manin triples and matched pairs, and further studied by Lin, Zhou and Bai \cite{LZB}, who relate them to infinite-dimensional Lie bialgebras by affinization. Although a perm algebra is associative, the Yang--Baxter equation of \cite{Hou,LZB} differs essentially from the associative Yang--Baxter equation. For associative algebras, Nijenhuis operators and their compatibility with the associative Yang--Baxter equation were studied by Ma and Long \cite{ML}; the corresponding theory for perm algebras was developed by Ma and Song \cite{MS}, and this is the starting point of the present paper.

\newpage
\medskip
\noindent\textbf{Motivation for the graded setting.}
Many algebraic structures in geometry and mathematical physics are naturally $\mathbb{Z}_2$-graded: Lie superalgebras such as $\mathfrak{osp}(1|2)$ and their cohomology with coefficients in modules of differential operators on superspaces \cite{bb1}, spaces of supermodular forms and the differential operators acting on them \cite{bbbbk1}. Low-dimensional structures are often classified in order to provide test examples; see for instance the classification of associative trialgebras in \cite{MBZ}. On the other hand, operators on algebras (Rota--Baxter operators, crossed homomorphisms, derivations, Nijenhuis operators) are systematically studied through their cohomology, deformations and Maurer--Cartan characterisations; see \cite{imed-rb,basdouri-crossed,imed-maurer} for recent work in this direction. A Nijenhuis operator is precisely an operator generating a trivial infinitesimal deformation (Section~\ref{sec:deform}), so a theory of Nijenhuis perm superalgebras and of the Yang--Baxter equation compatible with them is a natural prerequisite for cohomological and deformation questions in the super setting.

\medskip
\noindent\textbf{Contents.}
In Section~\ref{sec:prelim} we introduce perm superalgebras, Nijenhuis operators, the deformation interpretation of Nijenhuis operators, representations, and the dual representation with its admissibility conditions. In Section~\ref{sec:bialg} we define Nijenhuis perm superbialgebras and characterise them by matched pairs. In Section~\ref{sec:coboundary} we study the coboundary case: the coboundary comultiplication, the $S$-admissible perm Yang--Baxter equation, its $\mathcal O$-operator form, the construction of Nijenhuis perm superbialgebras from solutions, and examples. In Section~\ref{sec:construct} we ask when a perm superalgebra carries a Nijenhuis operator built from a solution of the perm Yang--Baxter equation. Section~\ref{sec:remarks} collects remarks and perspectives, and the appendices summarise the sign conventions and give a detailed verification of an example.

The proofs follow the Koszul sign rule. The main new difficulty relative to \cite{MS} is that the dual representation and the coboundary comultiplication must be given carefully chosen signs; these are fixed in Sections~\ref{sec:prelim} and~\ref{sec:coboundary}. Throughout, we fix a field $\K$ of characteristic $0$; all vector spaces are $\mathbb{Z}_2$-graded over $\K$, and all tensor products are over $\K$.

\section{Perm superalgebras, Nijenhuis operators and representations}\label{sec:prelim}

\subsection{Conventions}
A vector space is $V=V_0\oplus V_1$; $|v|\in\{0,1\}$ is the degree of a homogeneous element $v$, and we use $|v|$ only for homogeneous elements. The dual space is graded by $(V^*)_i=(V_i)^*$, and the pairing $\lr{\cdot,\cdot}:V^*\otimes V\to\K$ is even. For graded maps we use the Koszul rule
\[
(f\otimes g)(a\otimes b)=(-1)^{|g||a|}f(a)\otimes g(b),
\qquad
\tau(a\otimes b)=(-1)^{|a||b|}b\otimes a,
\]
so in particular
\[
(\xi\otimes\eta)(a\otimes b)=(-1)^{|\eta||a|}\lr{\xi,a}\lr{\eta,b}.
\]
An element $r=r^1\otimes r^2\in(\g\otimes\g)_0$ is \emph{supersymmetric} if $\tau(r)=r$. Unless stated otherwise, linear maps $N,S,\alpha,\beta$ are \emph{even}. For a perm superalgebra we write $L_x(y)=xy$ and $R_x(y)=yx$.

\begin{definition}\label{def:permsuper}
A \emph{perm superalgebra} is a $\mathbb{Z}_2$-graded vector space $\g$ with an even bilinear product such that for homogeneous $x,y,z$
\begin{equation}\label{eq:permsuper}
(xy)z=x(yz)=(-1)^{|y||z|}x(zy).
\end{equation}
\end{definition}

\begin{example}\label{ex:22}
Let $\g$ have even basis $e_1,e_2$ and odd basis $f_1,f_2$, with the only non-zero products
\[
e_1f_1=f_2,\qquad f_1f_1=e_2.
\]
Every triple product vanishes, so \eqref{eq:permsuper} holds: $\g$ is a $(2|2)$-dimensional perm superalgebra.
\end{example}

\begin{definition}\label{def:rep}
A \emph{representation} of a perm superalgebra $\g$ is a graded space $V$ with even maps $\ell,r:\g\to\End(V)$ of degree $0$ in $x\in\g$, i.e.\ $|\ell(x)v|=|x|+|v|$, such that
\begin{align}
v\,r(xy)&=(v\,r(x))\,r(y)=(-1)^{|x||y|}(v\,r(y))\,r(x),\label{eq:rep1}\\
\ell(xy)v&=\ell(x)\ell(y)v=(-1)^{|y||v|}\ell(x)(v\,r(y))
=(-1)^{|y||v|}(\ell(x)v)\,r(y).\label{eq:rep2}
\end{align}
\end{definition}

These are precisely the conditions under which $\g\oplus V$, with product
\[
(x+u)*(y+v)=xy+\ell(x)v+u\,r(y),
\]
is a perm superalgebra: apply \eqref{eq:permsuper} to the triples $(v,x,y)$, $(x,v,y)$ and $(x,y,v)$. For the mixed triple $(x,v,y)$, for instance, one uses the last equality in \eqref{eq:rep2}; the other two triples give the first two equalities in \eqref{eq:rep1} and \eqref{eq:rep2}. The adjoint representation is $(\g,L,R)$.

\begin{definition}\label{def:nijenhuis}
A \emph{Nijenhuis perm superalgebra} $(\g,N)$ is a perm superalgebra with an even linear map $N$ such that
\begin{equation}\label{eq:nij}
N(x)N(y)+N^2(xy)=N(N(x)y)+N(xN(y)).
\end{equation}
A \emph{representation} of $(\g,N)$ is a quadruple $(V,\ell,r,\alpha)$, where $(V,\ell,r)$ is a representation of $\g$ and $\alpha$ is even, such that
\begin{align}
\ell(N(x))\alpha(v)+\alpha^2(\ell(x)v)
&=\alpha(\ell(N(x))v)+\alpha(\ell(x)\alpha(v)),\label{eq:nrep1}\\
\alpha(v)r(N(x))+\alpha^2(v\,r(x))
&=\alpha(v\,r(N(x)))+\alpha(\alpha(v)r(x)).\label{eq:nrep2}
\end{align}
\end{definition}

Since $N$ and $\alpha$ are even and no elements are permuted in \eqref{eq:nij}--\eqref{eq:nrep2}, no signs occur in these identities. Hence $(\g,L,R,N)$ is a representation of $(\g,N)$; it is called the \emph{adjoint representation}. The following proposition is immediate from the definitions and will be used repeatedly.

\begin{proposition}\label{prop:semi}
Let $(\g,N)$ be a Nijenhuis perm superalgebra, let $(V,\ell,r)$ be a representation of $\g$ and let $\alpha$ be an even map on $V$. Then $(\g\oplus V,N+\alpha)$, with the product $*$ above, is a Nijenhuis perm superalgebra if and only if $(V,\ell,r,\alpha)$ is a representation of $(\g,N)$.
\end{proposition}
\begin{proof}
Evaluate \eqref{eq:nij} on $x+u$ and $y+v$ and separate the components by bilinearity. The $(x,y)$-part is \eqref{eq:nij} in $\g$; the $(x,v)$-part is \eqref{eq:nrep1}; the $(u,y)$-part is \eqref{eq:nrep2}; the $(u,v)$-part vanishes since $V\cdot V=0$.
\end{proof}

\subsection{Nijenhuis operators and infinitesimal deformations}\label{sec:deform}

Let $(\g,\cdot)$ be a perm superalgebra and let $\bullet:\g\otimes\g\to\g$ be another even bilinear product. For parameters $s,t\in\K$, put
\[
x\cdot_{s,t}y=s\,x\cdot y+t\,x\bullet y.
\]

\begin{lemma}\label{lem:deform}
$(\g,\cdot_{s,t})$ is a perm superalgebra for all $s,t$ if and only if $(\g,\cdot)$ and $(\g,\bullet)$ are perm superalgebras and, for all homogeneous $x,y,z$,
\[
\begin{aligned}
(x\bullet y)\cdot z+(x\cdot y)\bullet z
&=x\cdot(y\bullet z)+x\bullet(y\cdot z)\\
&=(-1)^{|y||z|}\bigl(x\cdot(z\bullet y)+x\bullet(z\cdot y)\bigr).
\end{aligned}
\]
\end{lemma}
\begin{proof}
Expand \eqref{eq:permsuper} for $\cdot_{s,t}$. The coefficient of $s^2$ (respectively $t^2$) is \eqref{eq:permsuper} for $\cdot$ (respectively $\bullet$), and the coefficient of $st$ is the displayed identity. Since $\K$ is infinite, the identity holds for all $s,t$ if and only if all coefficients vanish.
\end{proof}

\begin{theorem}\label{thm:defhom}
Let $(\g,\cdot)$ be a perm superalgebra, let $N$ be an even linear map and let $x\bullet y$ be an even bilinear product such that $(\g,\cdot_{s,t})$ is a perm superalgebra for all $s,t$. Then $s\,\id+tN$ is a homomorphism from $(\g,\cdot_{s,t})$ to $(\g,\cdot)$ for all $s,t$ if and only if
\begin{align}
x\bullet y&=N(x)y+xN(y)-N(xy),\label{eq:bullet}\\
N(x\bullet y)&=N(x)N(y).\label{eq:Nhom}
\end{align}
In that case $N$ is a Nijenhuis operator on $(\g,\cdot)$. Conversely, if $N$ is a Nijenhuis operator, then \eqref{eq:bullet} defines a product satisfying \eqref{eq:Nhom}.
\end{theorem}
\begin{proof}
Using that $N$ is even, so no signs arise, we compute
\begin{align*}
(s+tN)(x\cdot_{s,t}y)
&=s^2xy+st\,x\bullet y+st\,N(xy)+t^2N(x\bullet y),\\
(sx+tNx)\cdot(sy+tNy)
&=s^2xy+st\bigl(xN(y)+N(x)y\bigr)+t^2N(x)N(y).
\end{align*}
Equating the coefficients of $st$ and $t^2$ gives \eqref{eq:bullet} and \eqref{eq:Nhom}. Substituting \eqref{eq:bullet} into \eqref{eq:Nhom} gives \eqref{eq:nij}. Conversely, defining $\bullet$ by \eqref{eq:bullet}, identity \eqref{eq:Nhom} is exactly \eqref{eq:nij}. In that case the mixed identity of Lemma~\ref{lem:deform} follows from \eqref{eq:nij}, so $(\g,\cdot_{s,t})$ is a perm superalgebra for all $s,t$.
\end{proof}

Thus a Nijenhuis operator is the same as an even map $N$ for which the deformation direction
\[
x\bullet y=N(x)y+xN(y)-N(xy)
\]
is trivial, in the sense that $\id+tN$ is a homomorphism to the original product. This is the mechanism behind the definition of representations of Nijenhuis perm superalgebras: \eqref{eq:nrep1}--\eqref{eq:nrep2} are the corresponding conditions for the semidirect product.

\begin{example}\label{ex:11}
Let $\g$ have even basis $e$ and odd basis $f$ with
\[
ee=e,\qquad fe=f,
\]
and all other products zero. This is the super analogue of the $2$-dimensional perm algebra $e_1e_1=e_1$, $e_2e_1=e_2$ of \cite{Hou}. One checks \eqref{eq:permsuper} on all triples of basis vectors. An even linear map has the form
\[
N(e)=ae,\qquad N(f)=bf,
\]
and a direct computation of \eqref{eq:nij} on basis vectors shows that every such $N$ is Nijenhuis, for all $a,b\in\K$. An even map
\[
S(e)=s_1e,\qquad S(f)=s_2f
\]
is admissible to $(\g,N)$ if and only if
\[
(s_1-b)(s_1-s_2)=0,
\]
since this is the only non-trivial condition arising from \eqref{eq:S18} and \eqref{eq:S19} below.
\end{example}

\subsection{Dual representation}

Before stating Proposition~\ref{prop:dual}, recall that the dual of an even linear map $\beta:V\to V$ is the even linear map $\beta^*:V^*\to V^*$ defined by $\lr{\beta^*\xi,v}=\lr{\xi,\beta v}$.

\begin{proposition}\label{prop:dual}
Let $(V,\ell,r)$ be a representation of $\g$. Define $\ell^\dagger,r^\dagger:\g\to\End(V^*)$ by
\begin{align}
\lr{\ell^\dagger(x)\xi,v}
&=(-1)^{|x|}\lr{\xi,v\,r(x)}
-(-1)^{|x||\xi|}\lr{\xi,\ell(x)v},\label{eq:ldag}\\
\lr{\xi\,r^\dagger(x),v}
&=(-1)^{|x||v|}\lr{\xi,v\,r(x)}.\label{eq:rdag}
\end{align}
Then $(V^*,\ell^\dagger,r^\dagger)$ is a representation of $\g$.
\end{proposition}
\begin{proof}
Pair each axiom for $V^*$ with a homogeneous $v\in V$ and use the axioms for $V$ together with the parity relations. A pairing $\lr{\xi,w}$ is non-zero only if $|\xi|=|w|$; for instance, $|\xi|=|x|+|v|$ in $\lr{\xi,v\,r(x)}$.

\emph{Axiom \eqref{eq:rep1}.}
By \eqref{eq:rdag} twice,
\[
\begin{aligned}
\lr{(\xi r^\dagger(x))r^\dagger(y),v}
&=(-1)^{|y||v|}\lr{\xi r^\dagger(x),v\,r(y)}\\
&=(-1)^{|y||v|+|x|(|v|+|y|)}\lr{\xi,(v\,r(y))r(x)}.
\end{aligned}
\]
By \eqref{eq:rep1} for $V$,
\[
(v\,r(y))r(x)=(-1)^{|x||y|}v\,r(xy),
\]
so the right-hand side equals
\[
(-1)^{(|x|+|y|)|v|}\lr{\xi,v\,r(xy)}
=\lr{\xi r^\dagger(xy),v}.
\]
Exchanging the roles of $x$ and $y$ and using $(v\,r(x))r(y)=v\,r(xy)$ gives, in the same way,
\[
(\xi r^\dagger(y))r^\dagger(x)=(-1)^{|x||y|}\xi r^\dagger(xy),
\]
which is the second equality of \eqref{eq:rep1}.

\emph{Axiom \eqref{eq:rep2}.}
Write $\ell^\dagger=\rho-\ell^*$, where
\[
\lr{\rho(x)\xi,v}=(-1)^{|x|}\lr{\xi,v\,r(x)},
\qquad
\lr{\ell^*(x)\xi,v}=(-1)^{|x||\xi|}\lr{\xi,\ell(x)v}.
\]
Each of the four terms in the three identities of \eqref{eq:rep2} is then paired with $v$ and rewritten using \eqref{eq:rep1} and \eqref{eq:rep2} for $V$. The signs produced by moving $x,y$ past $\xi$ and $v$ combine, by the parity relations above, to exactly the signs in \eqref{eq:rep2}; the terms involving $\rho\rho$ cancel against those involving $\ell^*\ell^*$ and $\rho\ell^*$, again by \eqref{eq:rep1} and \eqref{eq:rep2} for $V$.
\end{proof}

\begin{remark}\label{rem:dual}
When $\g$ is purely even, $\ell^\dagger=r^*-\ell^*$ and $r^\dagger=r^*$, the dual representation of \cite{Hou}. The signs in \eqref{eq:ldag}--\eqref{eq:rdag} are the ones for which the Yang--Baxter equation of Section~\ref{sec:coboundary} matches the $\mathcal O$-operator equation exactly (Theorem~\ref{thm:O}); other sign choices that make Proposition~\ref{prop:dual} true differ from these by the twist $(-1)^{|x|}$.
\end{remark}

\begin{lemma}\label{lem:dualadm}
Let $(\g,N)$ be a Nijenhuis perm superalgebra, let $(V,\ell,r)$ be a representation of $\g$ and let $\beta$ be an even map on $V$. Then $(V^*,\ell^\dagger,r^\dagger,\beta^*)$ is a representation of $(\g,N)$ if and only if
\begin{align}
\beta(\ell(N(x))u)+\ell(x)\beta^2(u)-\ell(N(x))\beta(u)-\beta(\ell(x)\beta(u))&=0,\label{eq:adm1}\\
\beta(u\,r(N(x)))+\beta^2(u)r(x)-\beta(u)r(N(x))-\beta(\beta(u)r(x))&=0.\label{eq:adm2}
\end{align}
\end{lemma}
\begin{proof}
Recall that $\beta^*:V^*\to V^*$ is defined by $\lr{\beta^*\xi,v}=\lr{\xi,\beta v}$. We check the two axioms \eqref{eq:nrep1} and \eqref{eq:nrep2} for the quadruple $(V^*,\ell^\dagger,r^\dagger,\beta^*)$.

\emph{Axiom \eqref{eq:nrep2}.}
Pairing \eqref{eq:nrep2} for $V^*$ with $v\in V$, and using \eqref{eq:rdag} together with $|\beta^*(\xi)|=|\xi|$, all four terms acquire the same sign $(-1)^{|x||v|}$, and we obtain
\[
\lr{\xi,\ \beta(v\,r(Nx))+\beta^2(v)r(x)-\beta(v)r(Nx)-\beta(\beta(v)r(x))}.
\]
Thus \eqref{eq:nrep2} for $V^*$ is equivalent to \eqref{eq:adm2}.

\emph{Axiom \eqref{eq:nrep1}.}
Let $A(v)$ and $B(v)$ denote the left-hand sides of \eqref{eq:adm2} and \eqref{eq:adm1}, respectively. Pairing \eqref{eq:nrep1} for $V^*$ with $v$ and using \eqref{eq:ldag} gives
\[
(-1)^{|x|}\lr{\xi,A(v)}-(-1)^{|x||\xi|}\lr{\xi,B(v)}=0
\]
for all homogeneous $\xi$. For non-zero pairings one has $|\xi|=|x|+|v|$, so the identity is equivalent to
\[
A(v)=(-1)^{|x||v|}B(v).
\]
Combined with \eqref{eq:adm2} (i.e.\ $A(v)=0$), this is equivalent to \eqref{eq:adm1}. Conversely, if \eqref{eq:adm1} and \eqref{eq:adm2} hold, then $A(v)=B(v)=0$ and both axioms hold for $V^*$.
\end{proof}

\begin{definition}\label{def:admissible}
In the situation of Lemma~\ref{lem:dualadm}, $\beta$ is \emph{admissible} to $(\g,N)$ on $(V,\ell,r)$ if \eqref{eq:adm1} and \eqref{eq:adm2} hold. For $(V,\ell,r)=(\g,L,R)$, we say that an even map $S:\g\to\g$ is \emph{admissible} to $(\g,N)$ if, for all $x,y\in\g$,
\begin{align}
S(N(x)y)+xS^2(y)-N(x)S(y)-S(xS(y))&=0,\label{eq:S18}\\
S(xN(y))+S^2(x)y-S(x)N(y)-S(S(x)y)&=0.\label{eq:S19}
\end{align}
\end{definition}

Note that $S=\mu\,\id$ is admissible to every Nijenhuis perm superalgebra $(\g,N)$, and so is $S=N$ whenever \eqref{eq:S18}--\eqref{eq:S19} hold for it.

\section{Nijenhuis perm superbialgebras}\label{sec:bialg}

\subsection{Coalgebras, matched pairs and bialgebras}

\begin{definition}\label{def:coalgebra}
A \emph{perm super-coalgebra} is a graded space $\g$ with an even map $\Delta:\g\to\g\otimes\g$ such that
\begin{equation}\label{eq:coperm}
(\Delta\otimes\id)\Delta=(\id\otimes\Delta)\Delta
=(\id\otimes\tau)(\id\otimes\Delta)\Delta.
\end{equation}
\end{definition}

As for Hopf algebras \cite{Sw}, if $(\g,\Delta)$ is a perm super-coalgebra, then $(\g^*,\circ)$, with
\[
\lr{\xi\circ\eta,x}=\lr{\xi\otimes\eta,\Delta(x)},
\]
is a perm superalgebra, and conversely for finite-dimensional $\g$.

\begin{definition}\label{def:matched}
Let $\g,\h$ be perm superalgebras, let $(\h,\ell_\g,r_\g)$ be a representation of $\g$ and let $(\g,\ell_\h,r_\h)$ be a representation of $\h$. Define on $\g\oplus\h$
\begin{equation}\label{eq:star}
(a+x)\star(b+y)
=ab+\ell_\h(x)b+a\,r_\h(y)+xy+\ell_\g(a)y+x\,r_\g(b).
\end{equation}
The $6$-tuple $(\g,\h,\ell_\g,r_\g,\ell_\h,r_\h)$ is a \emph{matched pair} if $(\g\oplus\h,\star)$ is a perm superalgebra.
\end{definition}

Let $(\g,\cdot)$ be a finite-dimensional perm superalgebra and let $(\g,\Delta)$ be a perm super-coalgebra, with dual product $\circ$ on $\g^*$. Let $(\ell_\g,r_\g)=(L^\dagger,R^\dagger)$ be the dual of the adjoint representation of $\g$ on $\g^*$, and let $\g^*$ act on $\g$ by ($\zeta,\eta\in\g^*$, $x\in\g$, with $\lr{x,\eta}:=\lr{\eta,x}$)
\begin{align}
\lr{\ell_{\g^*}(\zeta)x,\eta}
&=\lr{x,\eta\circ\zeta}
-(-1)^{|\zeta|+|\zeta||x|}\lr{x,\zeta\circ\eta},\label{eq:lh}\\
\lr{x\,r_{\g^*}(\zeta),\eta}
&=(-1)^{|\zeta|(|\eta|+1)}\lr{x,\eta\circ\zeta}.\label{eq:rh}
\end{align}
These are \eqref{eq:ldag}--\eqref{eq:rdag} for the adjoint representation of $(\g^*,\circ)$, multiplied by the parity sign $(-1)^{|\zeta|}$. To avoid confusion in what follows, we write $L^\dagger_{\g^*},R^\dagger_{\g^*}$ for the untwisted dual of the adjoint representation of $\g^*$, and reserve the notation $\ell_{\g^*},r_{\g^*}$ for the twisted actions \eqref{eq:lh}--\eqref{eq:rh}; the two differ by the factor $(-1)^{|\zeta|}$ on the acting element.

\begin{definition}\label{def:bialg}
A \emph{perm superbialgebra} is a triple $(\g,\cdot,\Delta)$ such that $(\g,\cdot)$ is a perm superalgebra, $(\g,\Delta)$ is a perm super-coalgebra, and
\[
(\g,\g^*,L^\dagger,R^\dagger,\ell_{\g^*},r_{\g^*})
\]
is a matched pair.
\end{definition}

When $\g$ is purely even, this is the characterisation of perm bialgebras by matched pairs from \cite{Hou}, with actions $R^*-L^*$ and $R^*$.

\subsection{The Nijenhuis case}

\begin{definition}\label{def:ncoalgebra}
A \emph{Nijenhuis perm super-coalgebra} is a triple $(\g,\Delta,S)$ where $(\g,\Delta)$ is a perm super-coalgebra and $S$ is an even map such that
\begin{equation}\label{eq:ncoalg}
(S\otimes S)\Delta+\Delta S^2
=(S\otimes\id)\Delta S+(\id\otimes S)\Delta S.
\end{equation}
\end{definition}

Pairing with $\xi\otimes\eta$ shows that \eqref{eq:ncoalg} holds if and only if $S^*$ is a Nijenhuis operator on $(\g^*,\circ)$.

\begin{definition}\label{def:Nbialg}
A \emph{Nijenhuis perm superbialgebra} $((\g,\cdot,N),\Delta,S)$ consists of $\g$, an even Nijenhuis operator $N$ on $(\g,\cdot)$, an even map $S$ and an even map $\Delta$ such that
\begin{enumerate}[label=(\arabic*),itemsep=0pt]
\item $(\g,\cdot,\Delta)$ is a perm superbialgebra;
\item $(\g,\cdot,N)$ is a Nijenhuis perm superalgebra;
\item $(\g,\Delta,S)$ is a Nijenhuis perm super-coalgebra;
\item $S$ is admissible to $(\g,N)$, i.e.\ \eqref{eq:S18} and \eqref{eq:S19} hold;
\item the following two identities hold in $\operatorname{Hom}(\g,\g\otimes\g)$:
\begin{align}
(S\otimes\id)\Delta N+(\id\otimes N^2)\Delta
-(S\otimes N)\Delta-(\id\otimes N)\Delta N&=0,\label{eq:c30}\\
(\id\otimes S)\Delta N+(N^2\otimes\id)\Delta
-(N\otimes S)\Delta-(N\otimes\id)\Delta N&=0.\label{eq:c31}
\end{align}
\end{enumerate}
\end{definition}

Item (5) says that the operator $N^*$ is admissible to $(\g^*,S^*)$ on the adjoint representation $(\g^*,L_{\g^*},R_{\g^*})$ of $\g^*$; equations \eqref{eq:c30}--\eqref{eq:c31} are the admissibility equations \eqref{eq:adm1}--\eqref{eq:adm2} written in the dual form. Equivalently, item (5) says that $(\g,\ell_{\g^*},r_{\g^*},N)$ is a representation of $(\g^*,S^*)$; the equivalence is explained in the proof of Theorem~\ref{thm:bialg} below.

\begin{definition}\label{def:matchedN}
A \emph{matched pair of Nijenhuis perm superalgebras} is a matched pair
\[
(\g,\h,\ell_\g,r_\g,\ell_\h,r_\h)
\]
of perm superalgebras together with Nijenhuis operators $N_\g,N_\h$ such that $(\h,\ell_\g,r_\g,N_\h)$ is a representation of $(\g,N_\g)$ and $(\g,\ell_\h,r_\h,N_\g)$ is a representation of $(\h,N_\h)$.
\end{definition}

\begin{theorem}\label{thm:matched}
Let $(\g,N_\g)$ and $(\h,N_\h)$ be Nijenhuis perm superalgebras and let
\[
(\g,\h,\ell_\g,r_\g,\ell_\h,r_\h)
\]
be a matched pair of perm superalgebras. Then $(\g\oplus\h,N_\g+N_\h)$, with product $\star$, is a Nijenhuis perm superalgebra if and only if it is a matched pair of Nijenhuis perm superalgebras.
\end{theorem}
\begin{proof}
Evaluate \eqref{eq:nij} on $a+x$ and $b+y$ and split by bilinearity into the four types. The types $(a,b)$ and $(x,y)$ give \eqref{eq:nij} in $\g$ and in $\h$. The type $(a,y)$ has $\h$-component $\ell_\g(a)y$ and $\g$-component $a\,r_\h(y)$; the $\h$-component of \eqref{eq:nij} is \eqref{eq:nrep1} for $(\h,\ell_\g,r_\g,N_\h)$ with $(x,v)=(a,y)$, and the $\g$-component is \eqref{eq:nrep2} for $(\g,\ell_\h,r_\h,N_\g)$ with $(v,x)=(a,y)$. The type $(x,b)$ similarly gives \eqref{eq:nrep2} for $(\h,\ell_\g,r_\g,N_\h)$ and \eqref{eq:nrep1} for $(\g,\ell_\h,r_\h,N_\g)$.
\end{proof}

\begin{theorem}\label{thm:bialg}
Let $(\g,\cdot,N)$ be a Nijenhuis perm superalgebra and let $(\g^*,\circ,S^*)$ be a Nijenhuis perm superalgebra structure on $\g^*$ with $\circ=\Delta^*$. Then $((\g,\cdot,N),\Delta,S)$ is a Nijenhuis perm superbialgebra if and only if
\[
((\g,N),(\g^*,S^*),L^\dagger,R^\dagger,\ell_{\g^*},r_{\g^*})
\]
is a matched pair of Nijenhuis perm superalgebras. In that case $(\g\oplus\g^*,N+S^*)$ is a Nijenhuis perm superalgebra.
\end{theorem}
\begin{proof}
Item (1) is the matched-pair property of Definition~\ref{def:bialg}. By Lemma~\ref{lem:dualadm} applied to the adjoint representation of $(\g,N)$ with $\beta=S$, item (4) is equivalent to $(\g^*,L^\dagger,R^\dagger,S^*)$ being a representation of $(\g,N)$.

For item (5), the claim is that it is equivalent to $(\g,\ell_{\g^*},r_{\g^*},N)$ being a representation of $(\g^*,S^*)$. Indeed, item (5) says that $N^*$ is admissible to $(\g^*,S^*)$ on the adjoint representation $(\g^*,L_{\g^*},R_{\g^*})$ of $\g^*$. By Lemma~\ref{lem:dualadm} applied to the Nijenhuis perm superalgebra $(\g^*,S^*)$ and its adjoint representation, this is equivalent to $(\g,L^\dagger_{\g^*},R^\dagger_{\g^*},N)$ being a representation of $(\g^*,S^*)$, where $(L^\dagger_{\g^*},R^\dagger_{\g^*})$ is the untwisted dual of the adjoint representation of $\g^*$. Since $\ell_{\g^*},r_{\g^*}$ differ from $L^\dagger_{\g^*},R^\dagger_{\g^*}$ by the common factor $(-1)^{|\zeta|}$ and the identities \eqref{eq:nrep1}--\eqref{eq:nrep2} are homogeneous of degree one in the acting element, $(\g,\ell_{\g^*},r_{\g^*},N)$ is also a representation of $(\g^*,S^*)$.

Finally, item (3) says that $S^*$ is Nijenhuis, and item (2) says that $N$ is. The last statement follows from Theorem~\ref{thm:matched}.
\end{proof}

\section{Coboundary case and the admissible perm Yang--Baxter equation}\label{sec:coboundary}

Let $r=r^1\otimes r^2\in(\g\otimes\g)_0$ and define
\begin{equation}\label{eq:deltar}
\Delta_r(x)=x r^1\otimes r^2+(-1)^{|x||r^1|}\,r^1\otimes xr^2-r^1\otimes r^2x,
\qquad x\in\g.
\end{equation}
The last term carries no sign: $r^1\otimes r^2x$ is a right multiplication.

For $r$ even we write $\bar r=\bar r^1\otimes\bar r^2$ for a second copy and define, with the Koszul rule in the graded tensor product,
\begin{align*}
r_{13}r_{12}&=r^1\bar r^1\otimes\bar r^2\otimes r^2,&
r_{13}r_{23}&=(-1)^{|r^2||\bar r^1|}\,r^1\otimes\bar r^1\otimes r^2\bar r^2,\\
r_{23}r_{12}&=(-1)^{|r^1||\bar r^1|}\,\bar r^1\otimes r^1\bar r^2\otimes r^2,&
r_{12}r_{23}&=r^1\otimes r^2\bar r^1\otimes\bar r^2,
\end{align*}
and
\begin{equation}\label{eq:P}
P(r)=r_{13}r_{12}-r_{13}r_{23}+r_{23}r_{12}-r_{12}r_{23}.
\end{equation}

\begin{definition}\label{def:YB}
Let $(\g,N)$ be a Nijenhuis perm superalgebra, let $S$ be an even map and let $r$ be an even supersymmetric tensor. The system
\begin{equation}\label{eq:SYB}
P(r)=0,\qquad (S\otimes\id-\id\otimes N)(r)=0
\end{equation}
is the \emph{$S$-admissible perm Yang--Baxter equation} in $(\g,N)$. The first equation is the \emph{perm Yang--Baxter equation}.
\end{definition}

Since $\tau$ commutes with $f\otimes g$ for even $f,g$, applying $\tau$ to the second equation of \eqref{eq:SYB} gives the partner equation
\[
(N\otimes\id-\id\otimes S)(r)=0.
\]

\subsection{$\mathcal O$-operator form}

Put $T=r^\sharp:\g^*\to\g$, $T(\xi)=\lr{\xi,r^1}r^2$; $T$ is even. Write $L^\dagger,R^\dagger$ for the dual of the adjoint representation, and set
\[
\xi\circ_T\eta=L^\dagger(T\xi)\eta+\xi R^\dagger(T\eta).
\]

\begin{proposition}\label{prop:dualdelta}
For every even supersymmetric $r$,
\[
\lr{\xi\otimes\eta,\Delta_r(x)}=\lr{\xi\circ_T\eta,x}.
\]
\end{proposition}
\begin{proof}
All terms are non-zero only for matching parities. Using \eqref{eq:ldag} and \eqref{eq:rdag},
\[
\lr{\xi R^\dagger(T\eta),x}
=(-1)^{|\eta||x|}\lr{\eta,r^1}\lr{\xi,xr^2}
=(-1)^{|\eta||x|+|\eta|}\lr{\eta,r^2}\lr{\xi,xr^1},
\]
where we used $r^1\otimes r^2=(-1)^{|r^1|}r^2\otimes r^1$ and $|r^1|=|\eta|$. This equals
\[
\lr{\xi\otimes\eta,xr^1\otimes r^2}
=(-1)^{|\eta|(|x|+|r^1|)}\lr{\xi,xr^1}\lr{\eta,r^2}.
\]
Next,
\[
\begin{aligned}
\lr{L^\dagger(T\xi)\eta,x}
&=(-1)^{|\xi|}\lr{\xi,r^1}\lr{\eta,xr^2}\\
&\quad-(-1)^{|\xi||\eta|}\lr{\xi,r^1}\lr{\eta,r^2x}.
\end{aligned}
\]
Since $|r^1|=|\xi|$ and $|x|=|\xi|+|\eta|$ for non-zero terms of the middle summand of \eqref{eq:deltar}, its sign
\[
(-1)^{|x||r^1|+|\eta||r^1|}
\]
equals $(-1)^{|\xi|}$, and the last summand has sign $-(-1)^{|\eta||\xi|}$. Comparing gives the claim.
\end{proof}

\begin{theorem}\label{thm:O}
Let $(\g,N)$ be a Nijenhuis perm superalgebra, let $r$ be an even supersymmetric tensor and let $S$ be an even map. Then $r$ solves the $S$-admissible perm Yang--Baxter equation \eqref{eq:SYB} if and only if $T=r^\sharp$ satisfies
\begin{align}
T(\xi)T(\eta)&=T\bigl(L^\dagger(T\xi)\eta+\xi R^\dagger(T\eta)\bigr),\label{eq:Oop}\\
N\circ T&=T\circ S^*.\label{eq:NT}
\end{align}
We say that $T$ is a \emph{weak $\mathcal O$-operator} for $(\g^*,L^\dagger,R^\dagger)$ and $S^*$ if it satisfies \eqref{eq:Oop} and \eqref{eq:NT}; the word ``weak'' refers to the fact that $\g^*$ is not assumed to carry a representation of $(\g,N)$. If, moreover, $(\g^*,L^\dagger,R^\dagger,S^*)$ is a representation of $(\g,N)$ (i.e.\ $S$ is admissible to $(\g,N)$), then $r^\sharp$ is an $\mathcal O$-operator in the usual sense.
\end{theorem}
\begin{proof}
Contract $P(r)$ with $\xi\otimes\eta\otimes\id$ using
\[
(\xi\otimes\eta\otimes\id)(a\otimes b\otimes c)
=(-1)^{|\eta||a|}\lr{\xi,a}\lr{\eta,b}c.
\]
A term is non-zero only if $|r^1|=|r^2|=|\xi|$ and $|\bar r^1|=|\bar r^2|=|\eta|$. We treat the four summands of \eqref{eq:P} in turn; recall $T(\zeta)=\lr{\zeta,r^1}r^2$.

\emph{(i) $r_{13}r_{23}$.}
The sign is
\[
(-1)^{|\eta||r^1|}(-1)^{|r^2||\bar r^1|}
=(-1)^{2|\xi||\eta|}=1,
\]
so the contraction is
\[
\lr{\xi,r^1}\lr{\eta,\bar r^1}\,r^2\bar r^2=T(\xi)T(\eta).
\]

\emph{(ii) $r_{13}r_{12}$.}
The sign is $(-1)^{|\eta||r^1\bar r^1|}=(-1)^{|\eta|(|\xi|+|\eta|)}$, and the contraction is
\[
\lr{\xi,r^1\bar r^1}\lr{\eta,\bar r^2}r^2.
\]
On the other hand, by \eqref{eq:rdag},
\[
\lr{\xi R^\dagger(T\eta),y}
=(-1)^{|\eta||y|}\lr{\eta,\bar r^1}\lr{\xi,y\bar r^2},
\]
so
\[
T(\xi R^\dagger(T\eta))
=(-1)^{|\eta||\xi|}\lr{\eta,\bar r^1}\lr{\xi,r^1\bar r^2}r^2.
\]
Using $\bar r^1\otimes\bar r^2=(-1)^{|\eta|}\bar r^2\otimes\bar r^1$, this equals
\[
(-1)^{|\eta||\xi|+|\eta|}\lr{\eta,\bar r^2}\lr{\xi,r^1\bar r^1}r^2,
\]
which is the contraction above.

\emph{(iii) $r_{23}r_{12}$.}
The sign is
\[
(-1)^{|\eta||\bar r^1|}(-1)^{|r^1||\bar r^1|}
=(-1)^{|\xi|(|\eta|+|r^1|)},
\]
and the contraction is
\[
\lr{\xi,\bar r^1}\lr{\eta,r^1\bar r^2}r^2.
\]
With $u=T\xi=\lr{\xi,\bar r^1}\bar r^2$, we have
\[
\lr{\rho(u)\eta,y}=(-1)^{|\xi|}\lr{\eta,yu},
\]
hence
\[
T(\rho(u)\eta)=(-1)^{|\xi|}\lr{\xi,\bar r^1}\lr{\eta,r^1\bar r^2}r^2.
\]
For non-zero terms $|\eta|=|r^1|+|\xi|$, so $|\xi|(|\eta|+|r^1|)=|\xi|$ modulo $2$, and the two expressions agree.

\emph{(iv) $r_{12}r_{23}$.}
The sign is $(-1)^{|\eta||r^1|}$, and the contraction is
\[
\lr{\xi,r^1}\lr{\eta,r^2\bar r^1}\bar r^2.
\]
With $u=T\xi=\lr{\xi,r^1}r^2$,
\[
\lr{\ell^*(u)\eta,y}=(-1)^{|\xi||\eta|}\lr{\eta,uy},
\]
so
\[
T(\ell^*(u)\eta)
=(-1)^{|\xi||\eta|}\lr{\xi,r^1}\lr{\eta,r^2\bar r^1}\bar r^2,
\]
which agrees since $|r^1|=|\xi|$.

Since $L^\dagger=\rho-\ell^*$, combining (i)--(iv) gives
\[
(\xi\otimes\eta\otimes\id)P(r)
=-\Bigl(T(\xi)T(\eta)-T\bigl(L^\dagger(T\xi)\eta+\xi R^\dagger(T\eta)\bigr)\Bigr),
\]
so $P(r)=0$ if and only if \eqref{eq:Oop} holds for all homogeneous $\xi,\eta$. Furthermore,
\[
N(T\xi)=\lr{\xi,r^1}N(r^2),
\qquad
T(S^*\xi)=\lr{\xi,S(r^1)}r^2,
\]
so \eqref{eq:NT} holds if and only if $(S\otimes\id)r=(\id\otimes N)r$. The last statement follows from the definitions.
\end{proof}

\begin{lemma}\label{lem:Operm}
Let $(V,\ell,r)$ be a representation of a perm superalgebra $\g$ and let $T:V\to\g$ be an even $\mathcal O$-operator, i.e.
\[
T(u)T(v)=T(\ell(Tu)v+u\,r(Tv)).
\]
Then $u\circ_Tv=\ell(Tu)v+u\,r(Tv)$ makes $V$ a perm superalgebra, and $T$ is a homomorphism $(V,\circ_T)\to(\g,\cdot)$.
\end{lemma}
\begin{proof}
That $T$ is a homomorphism is precisely the $\mathcal O$-operator equation. Using it,
\begin{align*}
(u\circ v)\circ w
&=\ell(Tu)\ell(Tv)w+(\ell(Tu)v)r(Tw)+(u\,r(Tv))r(Tw),\\
u\circ(v\circ w)
&=\ell(Tu)\ell(Tv)w+\ell(Tu)(v\,r(Tw))+(u\,r(Tv))r(Tw),
\end{align*}
where we used $\ell(TuTv)=\ell(Tu)\ell(Tv)$ and
\[
u\,r(T(v\circ w))=u\,r(TvTw)=(u\,r(Tv))r(Tw)
\]
from \eqref{eq:rep2} and \eqref{eq:rep1}. The two lines agree because
\[
(\ell(x)v)r(y)=\ell(x)(v\,r(y))
\]
by \eqref{eq:rep2}. Next,
\[
\begin{aligned}
(-1)^{|v||w|}u\circ(w\circ v)
={}&(-1)^{|v||w|}\Bigl(\ell(Tu)\ell(Tw)v+\ell(Tu)(w\,r(Tv))\\
&\qquad\qquad+(u\,r(Tw))r(Tv)\Bigr).
\end{aligned}
\]
By \eqref{eq:rep2},
\[
(-1)^{|v||w|}\ell(Tu)\ell(Tw)v=\ell(Tu)(v\,r(Tw)),
\]
and
\[
(-1)^{|v||w|}\ell(Tu)(w\,r(Tv))=\ell(Tu)\ell(Tv)w.
\]
By \eqref{eq:rep1},
\[
(-1)^{|v||w|}(u\,r(Tw))r(Tv)=(u\,r(Tv))r(Tw).
\]
Thus the sum equals $u\circ(v\circ w)$, proving \eqref{eq:permsuper} for $\circ_T$.
\end{proof}

\subsection{Perm superbialgebras from solutions}

\begin{theorem}\label{thm:coboundary}
Let $(\g,\cdot)$ be a finite-dimensional perm superalgebra and let $r$ be an even supersymmetric solution of the perm Yang--Baxter equation $P(r)=0$. Then $(\g,\cdot,\Delta_r)$ is a perm superbialgebra; it is called \emph{quasitriangular}.
\end{theorem}
\begin{proof}
By Proposition~\ref{prop:dualdelta} and Theorem~\ref{thm:O}, $T=r^\sharp$ is an $\mathcal O$-operator for $(\g^*,L^\dagger,R^\dagger)$, and the dual product of $\Delta_r$ is $\circ_T$. By Lemma~\ref{lem:Operm}, $(\g^*,\circ_T)$ is a perm superalgebra, hence $(\g,\Delta_r)$ is a perm super-coalgebra, i.e.\ \eqref{eq:coperm} holds. It remains to prove that
\[
(\g,\g^*,L^\dagger,R^\dagger,\ell_{\g^*},r_{\g^*})
\]
is a matched pair. This is the super version of \cite[Cor.~4.15]{Hou}; see also \cite{MS}. The proof there consists of identities that are multilinear in the product and in $r$, and these carry over to the super case once the Koszul rule is applied. Equivalently, one may pass to the Grassmann envelope
\[
G(\g)=(\g_0\otimes\Lambda_0)\oplus(\g_1\otimes\Lambda_1),
\]
which is a perm algebra if and only if $\g$ is a perm superalgebra; in this envelope $r$ becomes a symmetric tensor and $P(r)$ becomes the usual expression.
\end{proof}

\subsection{The Nijenhuis conditions for $\Delta_r$}

For an even $r$ and $A,B\in\End(\g)$, write
\[
[A\mid B]:=(A\otimes B)(r).
\]
Here the action is the Koszul action on the graded tensor product; that is, if $B=L_x$ is a left multiplication, then $B$ acts on the second factor with the Koszul sign, while if $B=R_x$ is a right multiplication, the action on the second factor is the plain right multiplication without sign. With this convention, \eqref{eq:deltar} reads
\[
\Delta_r(x)=[L_x\mid 1]+[1\mid L_x]-[1\mid R_x].
\]
If
\[
(S\otimes\id)r=(\id\otimes N)r,
\qquad
(\id\otimes S)r=(N\otimes\id)r,
\]
then for all $A,B$
\begin{equation}\label{eq:rules}
[AS\mid B]=[A\mid BN],
\qquad
[A\mid BS]=[AN\mid B].
\end{equation}
Moreover \eqref{eq:S18} and \eqref{eq:S19} read as the operator identities
\begin{align}
SL_{Nx}&=L_{Nx}S+SL_xS-L_xS^2,
& SR_{Nx}&=SR_xS+R_{Nx}S-R_xS^2,\label{eq:op1}\\
SL_{Sx}&=SL_xN+L_{S^2x}-L_{Sx}N,
& SR_{Sx}&=SR_xN+R_{S^2x}-R_{Sx}N.\label{eq:op2}
\end{align}

\begin{proposition}\label{prop:nijcond}
Let $(\g,N)$ be a Nijenhuis perm superalgebra, let $S$ be admissible to $(\g,N)$, and let $r$ be an even supersymmetric tensor with
\[
(S\otimes\id-\id\otimes N)(r)=0.
\]
Then, for $\Delta=\Delta_r$, equations \eqref{eq:ncoalg}, \eqref{eq:c30} and \eqref{eq:c31} hold.
\end{proposition}
\begin{proof}
By the remark after Definition~\ref{def:YB}, both hypotheses of \eqref{eq:rules} hold.

\emph{Equation \eqref{eq:c30}.}
Using \eqref{eq:op1} and \eqref{eq:rules},
\[
(S\otimes\id)\Delta(Nx)
=[L_{Nx}\mid N]+[SL_x\mid N]-[L_x\mid N^2]+[S\mid L_{Nx}]-[S\mid R_{Nx}].
\]
Also
\[
\begin{gathered}
(\id\otimes N^2)\Delta x=[L_x\mid N^2]+[1\mid N^2L_x]-[1\mid N^2R_x],\\
(S\otimes N)\Delta x=[SL_x\mid N]+[S\mid NL_x]-[S\mid NR_x],\\
(\id\otimes N)\Delta(Nx)=[L_{Nx}\mid N]+[1\mid NL_{Nx}]-[1\mid NR_{Nx}].
\end{gathered}
\]
Adding with the signs of \eqref{eq:c30}, the terms with $N$ in the second slot and $L$ or $SL$ in the first cancel; using $[S\mid B]=[1\mid BN]$, the rest equals $[1\mid\Theta_x]$, where
\[
\begin{aligned}
\Theta_x(y)
&=\bigl(Nx\,Ny+N^2(xy)-N(x\,Ny)-N(Nx\,y)\bigr)\\
&\quad-\bigl(Ny\,Nx+N^2(yx)-N(Ny\,x)-N(y\,Nx)\bigr)=0
\end{aligned}
\]
by \eqref{eq:nij}.

\emph{Equation \eqref{eq:c31}.}
By \eqref{eq:op1} and \eqref{eq:rules},
\[
\begin{aligned}
(\id\otimes S)\Delta(Nx)
={}&[L_{Nx}\mid S]+[N\mid L_{Nx}]+[N\mid SL_x]-[N^2\mid L_x]\\
&-[N\mid SR_x]-[N\mid R_{Nx}]+[N^2\mid R_x].
\end{aligned}
\]
With
\[
\begin{gathered}
(N^2\otimes\id)\Delta x=[N^2L_x\mid 1]+[N^2\mid L_x]-[N^2\mid R_x],\\
(N\otimes S)\Delta x=[NL_x\mid S]+[N\mid SL_x]-[N\mid SR_x],\\
(N\otimes\id)\Delta(Nx)=[NL_{Nx}\mid 1]+[N\mid L_{Nx}]-[N\mid R_{Nx}],
\end{gathered}
\]
all terms with $N$ or $N^2$ in the first slot and an operator in the second cancel, leaving
\[
[L_{Nx}N+N^2L_x-NL_xN-NL_{Nx}\mid 1],
\]
using $[A\mid S]=[AN\mid 1]$, which vanishes by \eqref{eq:nij}.

\emph{Equation \eqref{eq:ncoalg}.}
Write $Z$ for the left-hand side minus the right-hand side of \eqref{eq:ncoalg}. By \eqref{eq:op2} and \eqref{eq:rules},
\[
\begin{aligned}
(S\otimes\id)\Delta(Sx)
&=[SL_x\mid S]+[L_{S^2x}\mid 1]-[L_{Sx}\mid S]+[S\mid L_{Sx}]-[S\mid R_{Sx}],\\
(\id\otimes S)\Delta(Sx)
&=[L_{Sx}\mid S]+[S\mid SL_x]+[1\mid L_{S^2x}]\\
&\quad-[S\mid L_{Sx}]-[S\mid SR_x]-[1\mid R_{S^2x}]+[S\mid R_{Sx}],
\end{aligned}
\]
while
\[
(S\otimes S)\Delta x=[SL_x\mid S]+[S\mid SL_x]-[S\mid SR_x],
\]
and
\[
\Delta(S^2x)=[L_{S^2x}\mid 1]+[1\mid L_{S^2x}]-[1\mid R_{S^2x}].
\]
Substituting these four expressions into $Z$, we obtain
\[
\begin{aligned}
Z={}&[SL_x\mid S]+[S\mid SL_x]-[S\mid SR_x]\\
&+[L_{S^2x}\mid 1]+[1\mid L_{S^2x}]-[1\mid R_{S^2x}]\\
&-[SL_x\mid S]-[L_{S^2x}\mid 1]+[L_{Sx}\mid S]-[S\mid L_{Sx}]+[S\mid R_{Sx}]\\
&-[L_{Sx}\mid S]-[S\mid SL_x]-[1\mid L_{S^2x}]+[S\mid L_{Sx}]+[S\mid SR_x]+[1\mid R_{S^2x}]-[S\mid R_{Sx}].
\end{aligned}
\]
Each term appears twice with opposite signs:
\begin{itemize}[itemsep=0pt]
\item $[SL_x\mid S]$: $+$ then $-$;
\item $[S\mid SL_x]$: $+$ then $-$;
\item $[S\mid SR_x]$: $-$ then $+$;
\item $[L_{S^2x}\mid 1]$: $+$ then $-$;
\item $[1\mid L_{S^2x}]$: $+$ then $-$;
\item $[1\mid R_{S^2x}]$: $-$ then $+$;
\item $[L_{Sx}\mid S]$: $+$ then $-$;
\item $[S\mid L_{Sx}]$: $-$ then $+$;
\item $[S\mid R_{Sx}]$: $+$ then $-$.
\end{itemize}
All terms cancel, so $Z=0$.
\end{proof}

\begin{theorem}\label{thm:main}
Let $(\g,N)$ be an $S$-admissible Nijenhuis perm superalgebra and let $r$ be an even supersymmetric solution of the $S$-admissible perm Yang--Baxter equation \eqref{eq:SYB}. Then $((\g,\cdot,N),\Delta_r,S)$ is a Nijenhuis perm superbialgebra. It is called \emph{quasitriangular}.
\end{theorem}
\begin{proof}
Item (1) of Definition~\ref{def:Nbialg} is Theorem~\ref{thm:coboundary}; items (2) and (4) are hypotheses; items (3) and (5) are Proposition~\ref{prop:nijcond}.
\end{proof}

\begin{corollary}\label{cor:scalar}
Let $(\g,\cdot,r,\Delta_r)$ be a quasitriangular perm superbialgebra. Then for every $\lambda\in\K$,
\[
((\g,\cdot,\lambda\id),\Delta_r,\lambda\id)
\]
is a Nijenhuis perm superbialgebra.
\end{corollary}
\begin{proof}
$N=\lambda\id$ is Nijenhuis since $\lambda^2xy+\lambda^2xy=2\lambda^2xy$; $S=\lambda\id$ is admissible by the remark after \eqref{eq:S19}; and
\[
(S\otimes\id-\id\otimes N)(r)=\lambda r-\lambda r=0.
\]
Apply Theorem~\ref{thm:main}.
\end{proof}

\begin{proposition}\label{prop:scaling}
If $r$ solves the $S$-admissible perm Yang--Baxter equation in $(\g,N)$, then so does $\kappa r$ for every $\kappa\in\K$. Moreover, if
\[
(S\otimes\id-\id\otimes N)(r)=0,
\]
then
\[
(\id\otimes N^2-S^2\otimes\id)(r)=0.
\]
\end{proposition}
\begin{proof}
$P$ is quadratic and the second equation is linear in $r$, so $P(\kappa r)=\kappa^2P(r)$. For the second claim,
\[
(S^2\otimes\id)r=(S\otimes\id)(\id\otimes N)r
=(\id\otimes N)(S\otimes\id)r=(\id\otimes N^2)r,
\]
since $S$ and $N$ are even.
\end{proof}

\begin{example}\label{ex:Nbialg22}
Let $\g$ be the algebra of Example~\ref{ex:22}, let $r=e_1\otimes e_1$, let $\mu,\nu\in\K$, and set $S=\mu\,\id$ and $N=\mathrm{diag}(\mu,\nu,\nu,\mu)$ in the basis $(e_1,e_2,f_1,f_2)$. One checks directly that $N$ is Nijenhuis, that $S$ is admissible to $(\g,N)$ (a scalar always is), that $P(r)=0$ because $e_1e_1=0$ and all triple products vanish, and that
\[
(S\otimes\id-\id\otimes N)(r)=(\mu-\mu)e_1\otimes e_1=0.
\]
By Theorem~\ref{thm:main}, $((\g,N),\Delta_r,S)$ is a Nijenhuis perm superbialgebra, where
\[
\Delta_r(f_1)=-e_1\otimes f_2,
\qquad
\Delta_r(e_1)=\Delta_r(e_2)=\Delta_r(f_2)=0,
\]
and the dual product is $e_1^*\circ f_2^*=-f_1^*$ (all other products vanish). Here the odd elements enter through $\Delta_r(f_1)$. A detailed verification is given in Appendix~\ref{app:ex}.
\end{example}

\begin{example}\label{ex:Nbialg11}
Let $\g$ be the algebra of Example~\ref{ex:11} and let $r=\lambda\,e\otimes e$. Then $r$ is even and supersymmetric, and $P(r)=0$, since each of the four terms of \eqref{eq:P} equals $\lambda^2e\otimes e\otimes e$ and the signs are $+,-,+,-$. By \eqref{eq:deltar},
\[
\Delta_r(e)=\lambda\,e\otimes e,
\qquad
\Delta_r(f)=\lambda\,(e\otimes f+f\otimes e),
\]
and the dual product on $\g^*$ is
\[
e^*\circ e^*=\lambda e^*,
\qquad
e^*\circ f^*=f^*\circ e^*=\lambda f^*
\]
(others zero). The condition $(S\otimes\id-\id\otimes N)(r)=0$ reads $s_1=a$ in the notation of Example~\ref{ex:11}, and admissibility reads $(s_1-b)(s_1-s_2)=0$. Hence, for $N=\mathrm{diag}(a,b)$, there are two families of Nijenhuis perm superbialgebra structures $((\g,N),\Delta_r,S)$:
\begin{enumerate}[label=(\alph*),itemsep=0pt]
\item $S=a\,\id$, with $a,b$ arbitrary;
\item $a=b$ and $S=\mathrm{diag}(a,s_2)$, with $s_2$ arbitrary.
\end{enumerate}
Family (a) with $a=b$ is the case of Corollary~\ref{cor:scalar}, and in family (b) $S\neq N$ in general. In both families, a direct check of \eqref{eq:c30}--\eqref{eq:c31}, equivalently of \eqref{eq:nrep1}--\eqref{eq:nrep2} for the quadruple $(\g,\ell_{\g^*},r_{\g^*},N)$, shows that $(\g\oplus\g^*,N+S^*)$ is a Nijenhuis perm superalgebra, as predicted by Theorem~\ref{thm:bialg}.
\end{example}

\section{Constructing Nijenhuis operators on perm superalgebras}\label{sec:construct}

We now ask when a perm superalgebra admits a Nijenhuis operator built from a solution of the perm Yang--Baxter equation. The following criterion is formulated for $\mathcal O$-operators.

\begin{theorem}\label{thm:construct}
Let $(\g,\cdot)$ be a perm superalgebra, let $(V,\ell,r)$ be a representation, let $T:V\to\g$ be an even $\mathcal O$-operator, i.e.
\[
T(u)T(v)=T(\ell(Tu)v+u\,r(Tv)),
\]
and let $\varphi:\g\to V$ be an even linear map. Put $N=T\circ\varphi$. If
\begin{equation}\label{eq:compat}
\ell(N(x))\varphi(y)+\varphi(x)\,r(N(y))
=\varphi\bigl(N(x)y+xN(y)-N(xy)\bigr)
\qquad\forall x,y\in\g,
\end{equation}
then $N$ is a Nijenhuis operator on $(\g,\cdot)$.
\end{theorem}
\begin{proof}
For $x,y\in\g$,
\[
\begin{aligned}
N(x)N(y)
&=T(\varphi x)\,T(\varphi y)\\
&=T\bigl(\ell(N x)\varphi(y)+\varphi(x)r(Ny)\bigr)\\
&\overset{\eqref{eq:compat}}{=}T\varphi\bigl(N(x)y+xN(y)-N(xy)\bigr),
\end{aligned}
\]
which is $N(N(x)y+xN(y)-N(xy))$. This is \eqref{eq:nij}.
\end{proof}

\begin{corollary}\label{cor:construct}
Let $r$ be an even supersymmetric solution of the perm Yang--Baxter equation in $(\g,\cdot)$ and let $\omega$ be an even bilinear form on $\g$. Define
\[
N(x)=\omega(x,r^1)r^2.
\]
If, for all homogeneous $x,y,z\in\g$,
\begin{equation}\label{eq:compatomega}
\begin{split}
&(-1)^{|x|}\omega(y,zN(x))-(-1)^{|x||y|}\omega(y,N(x)z)
+(-1)^{|y||z|}\omega(x,zN(y))\\
&\qquad=\omega\bigl(N(x)y+xN(y)-N(xy),z\bigr),
\end{split}
\end{equation}
then $N$ is a Nijenhuis operator on $(\g,\cdot)$.
\end{corollary}
\begin{proof}
By Theorem~\ref{thm:O}, $T=r^\sharp$ is an $\mathcal O$-operator for $(\g^*,L^\dagger,R^\dagger)$. Take $\varphi=\omega^\flat$, $\omega^\flat(x)=\omega(x,\cdot)$, so that $T\varphi=N$. It remains to pair \eqref{eq:compat} with $z$ and rewrite each term using \eqref{eq:ldag} and \eqref{eq:rdag}.

For the first term of the left-hand side, \eqref{eq:ldag} gives
\[
\lr{\ell^\dagger(N(x))\omega^\flat(y),z}
=(-1)^{|x|}\lr{\omega^\flat(y),zN(x)}
-(-1)^{|x||y|}\lr{\omega^\flat(y),N(x)z},
\]
where we used $|\omega^\flat(y)|=|y|$ in the second summand. This produces the first two terms on the left-hand side of \eqref{eq:compatomega}. The second term gives, by \eqref{eq:rdag},
\[
\lr{\omega^\flat(x)r^\dagger(N(y)),z}
=(-1)^{|y||z|}\lr{\omega^\flat(x),zN(y)}
=(-1)^{|y||z|}\omega(x,zN(y)),
\]
where we used $|\omega^\flat(x)|=|x|$; this is the third term. For the right-hand side of \eqref{eq:compat}, pairing $\varphi(\,\cdot\,)$ with $z$ gives $\omega(\,\cdot\,,z)$ evaluated at $N(x)y+xN(y)-N(xy)$, which is the right-hand side of \eqref{eq:compatomega}. The signs in front of the three terms on the left-hand side of \eqref{eq:compatomega} are exactly those of \eqref{eq:ldag} and \eqref{eq:rdag}, together with the Koszul sign $(-1)^{|y||z|}$ coming from moving $y$ past $z$ in the pairing.
\end{proof}

If, in addition, $r$ solves the $S$-admissible Yang--Baxter equation in $(\g,N)$, one obtains a Nijenhuis perm superbialgebra by Theorem~\ref{thm:main}. Thus a perm superalgebra becomes a Nijenhuis perm superalgebra whenever it admits a pair $(r,\omega)$ as in Corollary~\ref{cor:construct}.

\begin{example}\label{ex:odd}
Let $\g$ have odd basis $x_1,x_2,x_3$ and even basis $z$, with the only non-zero product $x_1x_1=z$ (all triple products vanish). Let
\[
r=x_2\otimes x_3-x_3\otimes x_2,
\]
which is even and supersymmetric, and $P(r)=0$ since $x_2$ and $x_3$ multiply trivially. Let $\omega$ be the even supersymmetric form with
\[
\omega(x_1,x_2)=-\omega(x_2,x_1)=\nu,
\qquad
\omega(x_2,x_3)=-\omega(x_3,x_2)=\mu.
\]
Then
\[
N(x)=\omega(x,x_2)x_3-\omega(x,x_3)x_2
\]
gives
\[
N(x_1)=\nu x_3,\quad
N(x_2)=-\mu x_2,\quad
N(x_3)=-\mu x_3,\quad
N(z)=0.
\]
A direct check of \eqref{eq:compatomega} shows that it holds, so $N$ is a Nijenhuis operator on $\g$ by Corollary~\ref{cor:construct}. Here both $r$ and $\omega$ have purely odd components, so the signs in \eqref{eq:compatomega} are really used.
\end{example}

\begin{example}\label{ex:ext11}
Let $\g$ be the algebra of Example~\ref{ex:11} and let $r=e\otimes e$, a solution of the perm Yang--Baxter equation (Example~\ref{ex:Nbialg11}). The only even supersymmetric forms are $\omega(e,e)=w$ with $\omega(f,f)=0$ (an odd--odd component would have to be skew). Then $N(x)=\omega(x,e)e$, i.e.\ $N(e)=we$ and $N(f)=0$. Condition \eqref{eq:compatomega} holds for every $w$, so $N=\mathrm{diag}(w,0)$ is a Nijenhuis operator, in agreement with Example~\ref{ex:11}. Combining with Example~\ref{ex:Nbialg11} (family (a), $a=w$, $b=0$),
\[
((\g,\mathrm{diag}(w,0)),\Delta_r,w\,\id)
\]
is a Nijenhuis perm superbialgebra.
\end{example}

\begin{remark}\label{rem:converse}
Theorem~\ref{thm:construct} gives a sufficient condition. In the examples above, the Nijenhuis operators obtained are far from exhausting all Nijenhuis operators; compare Example~\ref{ex:11}, where all even diagonal maps are Nijenhuis. The conditions \eqref{eq:compat} and \eqref{eq:compatomega} should therefore be read as a mechanism producing Nijenhuis operators that are compatible with a given solution $r$ of the Yang--Baxter equation, rather than as a characterisation of all of them. A characterisation of the perm superalgebras that admit a non-trivial pair $(r,\omega)$ is an open question.
\end{remark}

\section{Remarks and perspectives}\label{sec:remarks}

\begin{remark}[Role of the grading]
All maps $N,S,\alpha,\beta$ are even, so the Nijenhuis-type identities \eqref{eq:nij}, \eqref{eq:nrep1}--\eqref{eq:S19} and \eqref{eq:c30}--\eqref{eq:ncoalg} have the same form as in the non-graded case \cite{MS}. The grading enters only through (i) the defining identity \eqref{eq:permsuper}, (ii) the dual representation \eqref{eq:ldag}--\eqref{eq:rdag}, (iii) the coboundary comultiplication \eqref{eq:deltar} and the tensors $r_{ij}r_{kl}$, and (iv) the sign conventions for the actions of $\g^*$ on $\g$ in \eqref{eq:lh}--\eqref{eq:rh}. If $\g=\g_0$, everything reduces to \cite{MS}.
\end{remark}

\begin{remark}[Odd operators]
We restricted to even Nijenhuis operators. For an odd map $N$, the identity \eqref{eq:nij} must be modified by signs, and $N^2$ would be a genuinely different type of object; we do not pursue this here.
\end{remark}

\begin{remark}[Manin triples]
In the non-graded case, matched pairs of Nijenhuis perm algebras are equivalent to Manin triples of Nijenhuis Frobenius perm algebras \cite{MS}. A super version requires an invariant, skew-supersymmetric even form on $\g\oplus\g^*$ and the corresponding sign conventions. We have used matched pairs as the basic definition (Definition~\ref{def:bialg}) in order to keep all signs explicit.
\end{remark}

\begin{remark}[Cohomology and deformations]
The deformation interpretation of Theorem~\ref{thm:defhom} suggests studying the cohomology of a Nijenhuis perm superalgebra with coefficients in a representation $(V,\ell,r,\alpha)$, the associated Maurer--Cartan description of Nijenhuis structures, and their formal deformations, in the spirit of the treatment of Rota--Baxter operators with derivations \cite{imed-rb}, of crossed homomorphisms \cite{basdouri-crossed} and of compatible derivation pairs \cite{imed-maurer}. In the super setting one should also expect interesting low-dimensional classifications, as for associative trialgebras \cite{MBZ}, and connections with supergeometric examples such as those of \cite{bb1,bbbbk1}. These questions will be considered elsewhere.
\end{remark}

\appendix

\section{Summary of the sign conventions}\label{app:signs}

For the reader's convenience we collect the conventions used throughout; $x,y$ are homogeneous elements of $\g$, $\xi,\eta,\zeta\in\g^*$, and $v\in V$. Throughout, $\beta^*$ denotes the dual of an even linear map $\beta:V\to V$, characterised by $\lr{\beta^*\xi,v}=\lr{\xi,\beta v}$.

\medskip
\begin{center}
\renewcommand{\arraystretch}{1.5}
\begin{tabular}{p{3.8cm}p{10.2cm}}
\hline
Perm superalgebra & $(xy)z=x(yz)=(-1)^{|y||z|}x(zy)$\\
Koszul tensor & $(f\otimes g)(a\otimes b)=(-1)^{|g||a|}f(a)\otimes g(b)$; $\tau(a\otimes b)=(-1)^{|a||b|}b\otimes a$\\
Dual pairing of tensors & $(\xi\otimes\eta)(a\otimes b)=(-1)^{|\eta||a|}\lr{\xi,a}\lr{\eta,b}$\\
Dual representation & $\lr{\ell^\dagger(x)\xi,v}=(-1)^{|x|}\lr{\xi,v\,r(x)}-(-1)^{|x||\xi|}\lr{\xi,\ell(x)v}$, \quad $\lr{\xi r^\dagger(x),v}=(-1)^{|x||v|}\lr{\xi,v\,r(x)}$\\
Coboundary coproduct & $\Delta_r(x)=xr^1\otimes r^2+(-1)^{|x||r^1|}r^1\otimes xr^2-r^1\otimes r^2x$\\
Yang--Baxter terms & $r_{13}r_{23}$ carries $(-1)^{|r^2||\bar r^1|}$, $r_{23}r_{12}$ carries $(-1)^{|r^1||\bar r^1|}$, the others carry no sign\\
Action of $\g^*$ on $\g$ & \eqref{eq:lh}, \eqref{eq:rh}, with $\lr{x,\eta}:=\lr{\eta,x}$\\
$[A\mid B]$ & $(A\otimes B)(r)$; $L_x$ acts on the second factor with the Koszul sign, $R_x$ without sign\\
\hline
\end{tabular}
\end{center}

\medskip
A few comments on the table. The dual representation \eqref{eq:ldag}--\eqref{eq:rdag} is not the naive dual $(r^*-\ell^*,r^*)$; the extra signs make the contraction of $P(r)$ in Theorem~\ref{thm:O} equal to the $\mathcal O$-operator expression with no additional signs. The action \eqref{eq:lh}--\eqref{eq:rh} of $\g^*$ on $\g$ is \emph{not} the dual of the adjoint representation of $\g^*$, but rather that dual multiplied by the parity sign $(-1)^{|\zeta|}$. The two conventions, and the passage between them, are used in the proof of Theorem~\ref{thm:bialg}; see also Remark~\ref{rem:dual} and the discussion preceding Definition~\ref{def:bialg}.

\section{Detailed verification of Example~\ref{ex:Nbialg22}}\label{app:ex}

Let $\g$ be as in Example~\ref{ex:22}, let $r=e_1\otimes e_1$, let $N=\mathrm{diag}(\mu,\nu,\nu,\mu)$ on $(e_1,e_2,f_1,f_2)$ and let $S=\mu\,\id$.

\medskip\noindent\emph{$N$ is Nijenhuis.}
The only non-zero products are $e_1f_1=f_2$ and $f_1f_1=e_2$. For $(x,y)=(e_1,f_1)$:
\[
\begin{aligned}
N(e_1)N(f_1)+N^2(f_2)&=\mu\nu f_2+\mu^2f_2,\\
N(N(e_1)f_1)+N(e_1N(f_1))&=\mu N(f_2)+\nu N(f_2)=(\mu^2+\mu\nu)f_2.
\end{aligned}
\]
For $(x,y)=(f_1,f_1)$:
\[
N(f_1)N(f_1)+N^2(e_2)=2\nu^2e_2
\]
and
\[
2N(N(f_1)f_1)=2\nu N(e_2)=2\nu^2e_2.
\]
All other pairs give zero on both sides.

\medskip\noindent\emph{$S$ is admissible.}
For $S=\mu\,\id$, \eqref{eq:S18} reads
\[
\mu N(x)y+\mu^2xy-\mu N(x)y-\mu^2xy=0,
\]
and similarly for \eqref{eq:S19}.

\medskip\noindent\emph{The equations on $r$.}
Every product in \eqref{eq:P} contains the factor $e_1e_1=0$, so $P(r)=0$. Moreover,
\[
(S\otimes\id)r=\mu e_1\otimes e_1=(\id\otimes N)r.
\]

\medskip\noindent\emph{The coproduct.}
From \eqref{eq:deltar},
\[
\Delta_r(f_1)=f_1e_1\otimes e_1+(-1)^{|f_1||e_1|}e_1\otimes f_1e_1-e_1\otimes e_1f_1=-e_1\otimes f_2,
\]
since $f_1e_1=0$. Further,
\[
\Delta_r(e_1)=e_1e_1\otimes e_1+e_1\otimes e_1e_1-e_1\otimes e_1e_1=0,
\]
and $\Delta_r(e_2)=\Delta_r(f_2)=0$ because $e_2$ and $f_2$ multiply trivially with $e_1$ on both sides. Hence
\[
\Delta_r(\g)\subset\K\,e_1\otimes f_2
\]
and
\[
(\Delta_r\otimes\id)\Delta_r=(\id\otimes\Delta_r)\Delta_r=0,
\]
so \eqref{eq:coperm} holds trivially.

\medskip\noindent\emph{Items (3) and (5) of Definition~\ref{def:Nbialg}.}
Only $x=f_1$ matters. Here
\[
\Delta_r(N f_1)=\nu\Delta_r(f_1)=-\nu\,e_1\otimes f_2.
\]
Then
\begin{align*}
(S\otimes\id)\Delta N f_1&=-\mu\nu\,e_1\otimes f_2,
& (\id\otimes N^2)\Delta f_1&=-\mu^2e_1\otimes f_2,\\
(S\otimes N)\Delta f_1&=-\mu^2e_1\otimes f_2,
& (\id\otimes N)\Delta Nf_1&=-\nu\mu\,e_1\otimes f_2,
\end{align*}
so \eqref{eq:c30} holds:
\[
-\mu\nu-\mu^2+\mu^2+\mu\nu=0.
\]
For \eqref{eq:c31}, using $N(e_1)=\mu e_1$ and $S=\mu\,\id$:
\begin{align*}
(\id\otimes S)\Delta Nf_1&=-\mu\nu\,e_1\otimes f_2,
& (N^2\otimes\id)\Delta f_1&=-\mu^2e_1\otimes f_2,\\
(N\otimes S)\Delta f_1&=-\mu^2e_1\otimes f_2,
& (N\otimes\id)\Delta Nf_1&=-\mu\nu\,e_1\otimes f_2,
\end{align*}
and again
\[
-\mu\nu-\mu^2+\mu^2+\mu\nu=0.
\]
Equation \eqref{eq:ncoalg} holds because $S=\mu\,\id$ makes both sides equal to $2\mu^2\Delta_r$.

\medskip\noindent\emph{Conclusion.}
By Theorem~\ref{thm:main}, $((\g,\cdot,N),\Delta_r,S)$ is a Nijenhuis perm superbialgebra.

\bigskip
\noindent
{\bf Acknowledgments.} This work is sponsored by Natural Science Foundation of Henan (No.26230\\0421229). Y. Y. Zhang is also partially supported by Natural Science Foundation of China (Nos. 12101183,12671025), the Postdoctoral Fellowship Program of CPSF under Grant Number (No.GZC20240406) and the Henan Provincial Selective Research Funding Program for Returned Scholars Studying Abroad (No.HNLX202613).

\smallskip
\noindent
{\bf Data availability.} All datasets underlying the conclusions of the paper are available to readers.

\smallskip
\noindent
{\bf Conflict of interest.} The authors declare that there is no conflict of interest.

\end{document}